\documentclass[12pt,oneside,reqno]{amsart}
\usepackage{}
\usepackage{amssymb}
\usepackage{bbm}
\usepackage{cases}
\usepackage{amsmath}
\usepackage{graphicx}
\usepackage{mathrsfs}
\usepackage{stmaryrd}
\usepackage{amsfonts}
\usepackage{enumerate,amsmath,amssymb,amsthm}

\numberwithin{equation}{section}

\newcommand{\be}{\begin{eqnarray}}
\newcommand{\ee}{\end{eqnarray}}
\newcommand{\ce}{\begin{eqnarray*}}
\newcommand{\de}{\end{eqnarray*}}
\newtheorem{theorem}{Theorem}[section]
\newtheorem{lemma}[theorem]{Lemma}
\newtheorem{remark}[theorem]{Remark}
\newtheorem{definition}[theorem]{Definition}
\newtheorem{proposition}[theorem]{Proposition}
\newtheorem{Examples}[theorem]{Example}
\newtheorem{corollary}[theorem]{Corollary}

\def\p{\partial}

\def\[{{\Big[}}
\def\]{{\Big]}}
\def\<{{\langle}}
\def\>{{\rangle}}
\def\({{\Big(}}
\def\){{\Big)}}

\def\bx{{\mathbf{x}}}

\def\dif{{\mathord{{\rm d}}}}

\def\no{\nonumber}
\def\={&\!\!=\!\!&}
\def\bt{\begin{theorem}}
\def\et{\end{theorem}}
\def\bl{\begin{lemma}}
\def\el{\end{lemma}}
\def\br{\begin{remark}}
\def\er{\end{remark}}

\def\bd{\begin{definition}}
\def\ed{\end{definition}}
\def\bp{\begin{proposition}}
\def\ep{\end{proposition}}
\def\bc{\begin{corollary}}
\def\ec{\end{corollary}}
\def\bx{\begin{Examples}}
\def\ex{\end{Examples}}

\def\cI{{\mathcal I}}
\def\cJ{{\mathcal J}}

\def\cQ{{\mathcal Q}}

\def\mE{{\mathbb E}}

\def\mK{{\mathbb K}}

\def\mN{{\mathbb N}}

\def\mP{{\mathbb P}}

\def\mR{{\mathbb R}}

\def\sB{{\mathscr B}}

\def\sL{{\mathscr L}}

\def\sQ{{\mathscr Q}}

\def\geq{\geqslant}
\def\leq{\leqslant}

\allowdisplaybreaks

\begin{document}

\title{Heat kernel estimates for Dirichlet fractional Laplacian with gradient perturbation}

\date{}

\author{Peng Chen,\ \ Renming Song,\ \ Longjie Xie \ \  and \ \ Yingchao Xie}

\address{Peng Chen: Department of Mathematics, University of Macau,
Macau, P.R. China\\
Email: yb77430@umac.mo
 }

\address{Renming Song:
Department of Mathematics, University of Illinois,
Urbana, IL 61801, USA\\
Email: rsong@illinois.edu
 }

\address{Longjie Xie:
School of Mathematics and Statistics, Jiangsu Normal University,
Xuzhou, Jiangsu 221000, P.R.China\\
Email: xlj.98@whu.edu.cn
 }
\address{Yingchao Xie:
School of Mathematics and Statistics, Jiangsu Normal University,
Xuzhou, Jiangsu 221000, P.R.China\\
Email: ycxie@jsnu.edu.cn
 }

\thanks{Research of R. Song is supported by the Simons Foundation (\#429343, Renming Song). L. Xie is supported by NNSF of China (No. 11701233) and NSF of Jiangsu (No. BK20170226). Y. Xie is supported by NNSF of China (No. 11771187). The Project Funded by the PAPD of Jiangsu Higher Education Institutions is also gratefully acknowledged}

\begin{abstract}
We give a direct proof of the sharp two-sided estimates,
recently established in \cite{C-K-S-1, P-R}, for the Dirichlet heat
kernel of the fractional Laplacian with gradient perturbation in
$C^{1, 1}$ open sets  by using Duhamel formula.
We also obtain a gradient estimate for the Dirichlet heat kernel.
Our assumption on the open set is slightly weaker in that we only
require $D$ to be $C^{1,\theta}$  for some $\theta\in (\alpha/2, 1]$.
\bigskip


\end{abstract}

\maketitle \rm

\section{Introduction and main results}

Let $X=(X_t)_{t\geq 0}$ be an isotropic $\alpha$-stable
process on $\mathbb{R}^{d}$ with $d\geq1$ and $\alpha\in(0,2)$. The
infinitesimal generator of $X$ is the fractional Laplacian $\Delta^{\alpha/2}:=-(-\Delta)^{\alpha/2}$. For $f\in C^2_c(\mR^d)$, the fractional Laplacian $\Delta^{\alpha/2}$ can be written
in the following form:
$$
\Delta^{\alpha/2}f(x):= \int_{\mR^d}\big[f(x+z)-f(x)-1_{|z|\leq
1}z\cdot\nabla f(x)\big]\frac{c_{d,\alpha}}{|z|^{d+\alpha}}\dif z,
$$
where $c_{d,\alpha}$ is a positive constant.
It is well known
that the heat kernel $p(t,x,y)$ of $\Delta^{\alpha/2}$ (or
equivalently, the transition density of $X$) has the following
estimate:
\begin{align*}
p(t,x,y)\asymp \frac{t}{(|x-y|+t^{1/\alpha})^{d+\alpha}},\quad\forall (t,x,y)\in(0,\infty)\times\mR^d\times\mR^d.
\end{align*}
Here and below, for two non-negative functions $f$ and $g$, the notation $f \asymp g$ means that there are positive constants $c_1$ and $c_2$ such that $c_1g(x)\leq f(x)\leq c_2g(x)$ in the common domain of $f$ and $g$.

In \cite{Bo-Ja0}, Bogdan and Jakubowski studied the following perturbation of $\Delta^{\alpha/2}$ by a gradient operator
$$
\sL^b:=\Delta^{\alpha/2}+b\cdot\nabla
$$
in the case $d\ge2$ and $\alpha\in (1, 2)$.  They assumed that the  drift $b$  belongs to the Kato class defined below.

\bd
For any real-valued function $f$ on $\mR^d$, define for $r>0$
$$
K^\alpha_f(r):=\sup_{x\in
\mR^d}\int_{B(x,r)}\frac{|f(y)|}{|x-y|^{d+1-\alpha}}\dif y,
$$
where $B(x,r)$ denotes the open ball centered at
$x\in\mathbb{R}^{d}$ with radius $r$. Then $f$ is said to belong
to the Kato class $\mK^{\alpha-1}$ if $\lim_{r\downarrow
0}K^\alpha_f(r)=0$. \ed

In the remainder of this paper, we will
always assume that $d\geq2$ and $\alpha\in (1, 2)$, unless explicitly stated otherwise.
Intuitively, the
heat kernel $p^b(t,x,y)$ of $\sL^b$ should satisfy the following
Duhamel formula:
\begin{equation*}
p^b(t, x, y)=p(t, x, y)+\int^t_0\!\!\!\int_{\mR^d}p^b(t-s, x, z)b(z)\cdot\nabla_z p(s, z, y)\dif z
\dif s.
\end{equation*}
Define $p^b_0(t, x, y):=p(t, x, y)$ and for $k\geq 1$,
$$
p^b_k(t, x, y):=\int^t_0\!\!\!\int_{\mR^d}p^b_{k-1}(t-s, x,
z)b(z)\cdot\nabla_z p(s, z, y)\dif z\dif s.
$$
The following theorem is the main result of \cite{Bo-Ja0}.

\bt\label{t:BJ}
Assume that $b\in\mK^{\alpha-1}$.
\begin{enumerate}
\item There exist $T_0 > 0$ and $C > 1$ depending on $b$ only through the rate at which
$K^\alpha_{|b|} (r)$ goes to zero such that  $\sum^\infty_{k=0}p^b_k(t, x, y)$ converges locally
uniformly on $(0,T_0]\times \mR^d \times \mR^d$ to a positive jointly continuous function $p^b(t,x,y)$
and that on $(0,T_0]\times \mR^d \times \mR^d$,
$$
C^{-1}p(t, x, y)\leq p^b(t, x, y)\leq Cp(t, x, y).
$$
Moreover, $\int_{\mR^d} p^b(t,x,y)\dif y=1$ for every $t\in (0,T_0]$ and $x\in \mR^d$.

\item The function $p^b(t, x, y)$ can be extended uniquely to a positive
jointly continuous function on $(0, \infty) \times \mR^d \times \mR^d$ so that for all
$s, t \in  (0, \infty)$ and $(x, y)\in \mR^d\times \mR^d$, $\int_{\mR^d}  p^b(t,x,y)dy=1$
and
$$
p^b(t+s, x, y)=\int_{\mR^d}p^b(t, x, z)p^b(s, z, y)\dif z.
$$

\item If we define
$$
P^b_tf(x):=\int_{\mR^d}p^b(t, x, y)f(y)\dif y,
$$
then for any $f, g\in C^\infty_c(\mR^d)$,
$$
\lim_{t\downarrow0}\int_{\mR^d}t^{-1}\left(P^b_tf(x)-f(x) \right)g(x)\dif x=
\int_{\mR^d}(\sL^bf)(x)g(x)\dif x.
$$
Thus, $p^b(t, x, y)$ is the fundamental solution of $\sL^b$ in the
distributional sense.
\end{enumerate}
\et

Using the semigroup property, one can easily check that for any
$T>0$, there exists a constant $C> 1$ such that for all  $(t, x, y)\in(0, T]\times\mR^d\times\mR^d$,
\begin{equation}\label{e:newheat}
C^{-1} p(t, x, y)\le p^b(t, x, y)\le C p(t, x, y).
\end{equation}
It follows from \cite[Lemma 2.3]{C-K-S-1} that
$\{P^b_t, t\geq 0\}$ form a Feller semigroup, so there is a conservative Feller
process $X^b:=\{X^b_t, t\geq 0, \mP_x, x\in \mR^d\}$ on $\mR^d$ such
that $P^b_tf(x)=\mE_x[f(X^b_t)]$. The process $X^b$
is nonsymmetric and is called an $\alpha$-stable process with drift $b$.
See also \cite{C-H-X-Z,Xi-Zh} for the two-sided heat kernel estimates of more general non-local operators in the whole space $\mR^d$.

\vspace{2mm}
For any open subset $D \subset \mR^d$, define $\tau^b_D:= \inf\{t >
0:X^b_t\notin D\}$. We will use $X^{b,D}$ to denote the subprocess
of $X^b$ in $D$; that is, $X^{b,D}(\omega):=X^b(\omega)$ if $t
<\tau^b_D (\omega)$ and $X^{b,D}(\omega):=\partial$ if $t \geq
\tau^b_D(\omega)$, where $\partial$ is a cemetery state. Throughout
this paper, we use the convention that for every function $f$, we
extend its definition to $\partial$ by setting $f(\partial)=0$. The
infinitesimal generator of $X^{b,D}$ is given by $\sL^{b,
D}:=\sL^b|_D$, that is, $\sL^b$ on $D$ with zero exterior condition.
The processes $X^{b,D}$ has a joint continuous transition density
$p^{b, D} (t, x, y)$ which is also the Dirichlet heat kernel for
$\sL^{b, D}$. The subprocess of $X$ in $D$ will be denoted by $X^D$
and it is known to have a transition density $p^{D} (t, x, y)$.

Due to the complication near the boundary, sharp two-sided
estimates for the Dirichlet heat kernel are much more difficult to obtain.
The first sharp two-sided
estimates for the Dirichlet heat kernels of discontinuous Markov processes are due
to \cite{C-K-S-2}.
To state the related results, we first recall the definition of
$C^{1,\theta}$ open sets.
For $\theta\in(0,1]$,  an open set $D$ in
$\mR^{d}$ is said to be $C^{1,\theta}$ if there exist $r_{0}>0$ and
$\Lambda>0$ such that for every $Q\in\partial D$, there exist a
$C^{1,\theta}$-function $\phi=\phi_{Q}:
\mathbb{R}^{d-1}\rightarrow\mathbb{R}$ satisfying
$\phi(0)=\nabla\phi(0)=0$, $\|\nabla\phi\|_{\infty}\leq\Lambda$,
$|\nabla\phi(x)-\nabla\phi(z)|\leq\Lambda|x-z|^\theta$ and an
orthonormal coordinate system $y= (y_{1}, \cdots, y_{d-1},
y_{d}):=(\tilde{y}, y_{d})$ such that $B(Q, r_{0})\cap D=B(Q,
r_{0})\cap \{y: y_{d}>\phi(\tilde{y})\}$. The pair $(r_{0},
\Lambda)$ is called the characteristics of the $C^{1,\theta}$ open
set $D$. For $t>0$ and $x,y\in D$, we define
\begin{equation*}
q^D(t, x, y):=\Big(1\wedge\frac{\rho(x)^{\alpha/2}}{\sqrt{t}}\Big)\Big(1\wedge\frac{\rho(y)^{\alpha/2}}{\sqrt{t}}\Big)p(t,x,y),
\end{equation*}
where $\rho(x)$ denotes the distance between $x$ and $D^{c}$. In \cite{C-K-S-2},
Chen, Kim and Song proved that for any $d\geq 1$, $\alpha\in (0, 2)$ and
 $T>0$, when $D$ is a $C^{1,1}$ open set in $\mR^d$, there exists a constant $C>0$ such that
\begin{equation}\label{Heat2}
C^{-1}q^D(t, x, y)\leq p^D(t, x, y)\leq Cq^D(t, x, y), \qquad (t, x, y)\in (0, T]\times \mR^d\times
\mR^d.
\end{equation}
The above result has been generalized to $C^{1,\theta}$ open sets with
$\theta\in (\alpha/2, 1]$ in \cite{KK}.
As for the estimates of $p^{b, D} (t, x, y)$, the following result
is proved in \cite{C-K-S-1} in the case when $D$ is a bounded $C^{1,1}$ open set. The unbounded
case is due to \cite{P-R}.

\bt\label{t:CKS-KS}
Let $b\in\mK^{\alpha-1}$ and $D$ be a  $C^{1, 1}$ open set in $\mR^d$ with $C^{1, 1}$ characteristics
$(r_{0}, \Lambda)$. Then for any $T>0$, there exists a constant $C=C(T, r_0, \Lambda, d, \alpha, b)>1$ with the dependence on $b$ only via the rate at which $K^\alpha_{|b|}(r)$ tends to zero such that for all $(t, x, y)\in (0, T]\times \mR^d\times \mR^d$,
$$
C^{-1}q^D(t, x, y)\leq p^{b, D}(t, x, y)\leq Cq^D(t, x, y).
$$
\et

One might think that the estimates in Theorem \ref{t:CKS-KS} can be obtained from
the estimates \eqref{Heat2} for $p^D(t, x, y)$ using the following Duhamel formula:
\begin{align}\label{du}
p^{b,D}(t,x,y)=p^D(t,x,y)+\int_0^t\!\!\!\int_Dp^{b,D}(t-s,x,z)b(z)\cdot\nabla_z
p^D(s,z,y)\dif z\dif s.
\end{align}
However, unlike the whole space case, there was no good estimates on
$\nabla_z p^D(t, z, y)$, so the approach mentioned above could not
be carried through. Another obstacle to carrying out the approach
above in the present case is that
the following form of 3-P inequality (see \cite[Remark 2.3]{C-K-S-3}):
there exists $C>0$ such
that for any $0<s<t$ and $x,y,z\in D$,
\begin{align}
\frac{p^{D}(t-s,x,z)p^{D}(s,z,y)}{p^{D}(t,x,y)}\le C\big(
p^{D}(t-s,x,z)+p^{D}(s,z,y)\big),  \label{class3p}
\end{align}
does not hold (see \cite[Remark 2.3]{C-K-S-3}).
A whole space analog of the inequality above played a crucial role in proving the
estimates in Theorem \ref{t:BJ}. Partly due to the two reasons mentioned above, Theorem
\ref{t:CKS-KS} was much more difficult to prove than Theorem \ref{t:BJ}. To get around the
difficulties mentioned above,  \cite{C-K-S-1, P-R} used the Duhamel formula
for the Green functions
of $X^{b, D}$ and the probabilistic road-map designed in \cite{C-K-S-2} for establishing
the estimates \eqref{Heat2}.

\vspace{2mm}
In the recent paper \cite{K-R}, Kulczycki and Ryznar proved the following gradient estimate
for $p^D(t, x, y)$ (see \cite[Theorem 1.1 and Corollary 1.2]{K-R}): there exists a constant $C_1=C_1(d, \alpha)>0$ such that for any open set $D\subset \mR^d$ and all
$(t, x, y)\in (0, 1]\times D\times D$,
\begin{align*}
|\nabla_x p^{D}(t,x,y)|\leq  \frac{C_1}{\rho(x)\wedge
t^{1/\alpha}}p^{D}(t,x,y).
\end{align*}
It follows immediately that for any $T>0$,
there exists a constant $C_2=C_2(d, \alpha, T)>0$ such that for any open set $D\subset \mR^d$ and
all $(t, x, y)\in (0, T]\times D\times D$,
\begin{align}\label{grad2-a}
|\nabla_x p^{D}(t,x,y)|\leq  \frac{C_2}{\rho(x)\wedge
t^{1/\alpha}}p^{D}(t,x,y).
\end{align}
In this paper, we will use (\ref{grad2-a}) and the Duhamel formula \eqref{du}
to give a direct proof of Theorem \ref{t:CKS-KS}.
In fact, we will
establish two-sided estimates for $p^{b,D}$ with $b$ in a certain local Kato class and
$D$ being a $C^{1,\theta}$ open set with $\theta\in (\alpha/2, 1]$
instead of $C^{1,1}$ open set.
We also prove a gradient estimate for $p^{b,D}(t,x,y)$,
which is of independent interest.

\vspace{2mm}
To state our main results, we first introduce the following local
Kato class.

\bd
Let $D$ be any open set in $\mR^d$.
For any real-valued function $f$ defined on $D$, we define for every $r>0$,
$$
K^{\alpha, D}_f(r):=
\sup_{x\in D}\int_{D\cap
B(x,r)}\frac{|f(y)|}{|x-y|^{d+1-\alpha}} \dif y.
$$
Then $f$ is said to belong to the  local Kato class $\mK^{\alpha-1}_D$ if
$\lim_{r\downarrow0}K^{\alpha, D}_f(r)=0$.
\ed

\begin{remark}\label{r:Kato}
Using the maximum principle (see \cite[Theorem 5.2.2]{Ch}) it is easy to check that
a function $b: D\to \mR^d$ belongs to $\mK^{\alpha-1}_D$ if and only $b1_D$ belongs to
$\mK^{\alpha-1}$.
\end{remark}

The following is the main result of this paper.

\bt\label{main}
Let $D$ be a $C^{1,\theta}$ open subset of $\mathbb{R}^{d}$ with
$\theta\in(\alpha/2,1]$ and $b: D\to \mR^d$ belongs to $\mK^{\alpha-1}_D$.
Then there exists a unique function $p^{b,D}(t,x,y)$ on $(0,
\infty)\times \mR^d\times \mR^d$ satisfying \eqref{du} such that:
\begin{enumerate}[(i)]
\item for any $T>0$, there exists a constant $C_1>1$ such that for all $t\in(0,T]$ and $x, y\in D$, we have
\begin{align}\label{1.5}
C_1^{-1}q^D(t, x, y)\leq p^{b, D}(t, x, y)\leq C_1 q^D(t, x, y);
\end{align}

\item for any $T>0$, there exists a constant $C_2>0$ such that for all $t\in(0,T]$ and $x, y\in D$,
\begin{align}\label{vi}
|\nabla_x p^{b,D}(t,x,y)|\leq
\frac{C_2}{\rho(x)\wedge
t^{1/\alpha}}p^{D}(t,x,y),
\end{align}
and $p^{b,D}(t,x,y)$ also satisfies
\begin{align}\label{Duhamel}
p^{b,D}(t,x,y)=p^{D}(t,x,y)+\int_0^t\!\!\!\int_D
p^{D}(t-s,x,z)b(z)\cdot\nabla_z p^{b,D}(s,z,y)\dif
z\dif s;
\end{align}

\item for all $0<s<t$ and $x,y\in D$, the following Chapman-Kolmogorov's equation holds:
\begin{align}\label{eq21}
\int_{D}p^{b,D}(t-s,x,z)p^{b,D}(s,z,y)\dif z=p^{b,D}(t,x,y);
\end{align}

\item for any $f\in C^2_c(D)$, we have
\begin{align}\label{eq23}
P_{t}^{b,D}f(x)=f(x)+\int_0^tP_{t-s}^{b,D}\sL^{b,D}f(x)\dif s,
\end{align}
where $P_{t}^{b,D}f(x):=\int_{D}p^{b,D}(t,x,y)f(y)\dif y$;

\item for any $t\ge0$ and $x\in D$, it holds that
\begin{align}
\int_{D}p^{b,D}(t,x,y)\dif y\leq 1; \label{conser}
\end{align}

\item for any uniformly continuous function $f(x)$ with compact supports, we have
\begin{align}
\lim_{t\downarrow0} \|P_t^{b,D}f-f\|_\infty=0.   \label{con}
\end{align}
\end{enumerate}
\et

As an application of our heat kernel estimates,
we can get the following Harnack inequality on the semigroup $P_t^{b,D}$, which may be used to
study the long time behavior of the process, see, for eaxmple, \cite{L-L-W,W}.

\bc\label{cor}
There exists a constant $C$ such that for any non-negative function $f\in \sB_b(D)$, $T>0$ and $x,y\in D$, we have
\begin{align}
P_T^{b,D}f(x)\leq C\left(1\vee\frac{\rho(x)}{\rho(y)}\right)^{\alpha/2}\left(1\vee\frac{|x-y|}{(T\wedge1)^{1/\alpha}}\right)^{d+\alpha}P_T^{b,D}f(y).    \label{ha}
\end{align}
\ec

The remainder of this paper is organized as follows. In Section 2, we
prepare some important inequalities for latter use; the proof of
main result, Theorem \ref{main}, will be given in Section 3.

We conclude this introduction by spelling out some conventions that
will be used throughout this paper. The letter C with or without
subscripts will denote an unimportant constant and $f\preceq g$
means that $f\leq Cg$ for some $C\geq1$. The letter $\mathbb{N}$
will denote the collection of positive integers, and
$\mathbb{N}_{0}=\mathbb{N}\cup\{0\}$. We will use $:=$ to denote a
definition, $\sB_b(D)$ to denote the space of all bounded Borel measurable
functions on $D$ and we assume that all the functions
considered in this paper are Borel measurable.

\section{Preliminaries}

By combining \cite[Corollary 12]{Bo-Ja0} with Remark \ref{r:Kato}, we immediately get
the following equivalent characterization of $\mK^{\alpha-1}_D$,
which will be used in the proof of our heat kernel estimates.

\bl\label{include}
Let
$\beta>\frac{\alpha-1}{\alpha}.$ A function f belongs to ${\mK}^{\alpha-1}_D$ if and only if
\begin{equation*}
\lim_{t\rightarrow0}\sup_{x\in D}\int_{D}\bigg(\frac{1}{|y-x|^{d+1-\alpha}}\wedge\frac{t^{\beta}}{|y-x|^{d+1-\alpha+\alpha\beta}}\bigg)|f(y)|\dif y=0.
\end{equation*}
\el

The next result says that if $b\in\mK^{\alpha-1}$, then the density $p^{b,D}(t,x,y)$ of $X_t^{b,D}$ do satisfy the Duhamel formula \eqref{du}.
This result will used in the proof of our main result.

\bl\label{dens}
Assume that $b\in\mK^{\alpha-1}$ and $D$ is an open subset of $\mR^d$.
Then the transition density $p^{b,D}(t,x,y)$ of $X_t^{b,D}$ satisfies \eqref{du}.
\el

\begin{proof}
Let $\phi\in C_{c}^{\infty}\left((0,\infty)\times\mathbb{R}^{d}\right)$ with  ${\rm Supp}[\phi]\subset(0,1)\times B(0,1)$ and $\int_{0}^{\infty}\!\!\int_{\mathbb{R}^{d}}\phi(r,y)\dif y\dif r=1$. Fix $t>0$, for any $\psi\in C_{c}(D)$, define $f(s,x):=P^{D}_{t-s}\psi(x)$ and $f_{n}:=\phi_{n}\ast f$, where $\phi_{n}(r,y)=n^{d+1}\phi(nr, ny)$.
Let $D_j$ be a sequence of relatively compact open subsets of $D$ such that
$D_j\subset\overline{D_j}\subset D_{j+1}$ for all $j\ge 1$ and $D_j\uparrow D$.
Let $\tau_{D_j}^b:=\inf\{t>0: X_{t}^{b}\in D^{c}_j\}$.
It follows from \cite{C-W} that $X^b$ is a weak solution of the the stochastic differential equation
$$
dX^b_t=dX_t+b(X^b_t)dt.
$$
Thus by It\^o's formula, we have for sufficiently large $n$,
\begin{align*}
\mE\big[f_{n}\big(t\wedge\tau_{D_{j}}^{b},X_{t\wedge\tau_{D_{j}}^{b}}\big)\big]-f_{n}(0,x)
&=\int_{0}^{t}P_{s}^{b,D_{j}}\big[\p_s f_{n}+\Delta^{\alpha/2}f_{n}+b\cdot\nabla f_{n}\big](s,x)\dif s\\
&=\int_{0}^{t}P_{s}^{b,D_{j}}\big[\phi_{n}\ast\big(\p_s f+\Delta^{\alpha/2}f)+b\cdot\phi_{n}\ast\nabla f\big](s,x)\dif s.
\end{align*}

Since $b\in \mK^{\alpha-1}$ and
$$
p^{b, D_j}(t, x, y)\le p^b(t, x, y),
$$
applying \eqref{grad2-a}, \eqref{e:newheat} and
letting $n\rightarrow\infty$ we get
\begin{equation*}
\mE\big[f\big(t\wedge\tau_{D_{j}}^{b},X_{t\wedge\tau_{D_{j}}^{b}}\big)\big]=f(0,x)+\int_{0}^{t}P_{s}^{b,D_{j}}\big(b\cdot\nabla f\big)(s,x)\dif s.
\end{equation*}
Note that $f(t,x)=\psi(x)$ and $f(0,x)=P_{t}^{D}\psi(x),$ taking $j\rightarrow\infty$ we arrive at
$$
P_{t}^{b,D}\psi(x)=P_{t}^{D}\psi(x)+\int_{0}^{t}P_{s}^{b,D}\big(b\cdot\nabla P_{t-s}^{D}\psi\big)(x)\dif s,
$$
which in turn means the desired result.
\end{proof}

For any $\gamma\in\mR$,  we define for $t>0$ and $x\in \mR^d$,
$$
\varrho^\gamma(t,x):=\frac{t^\gamma}{(|x|+t^{1/\alpha})^{d+\alpha}}
$$
and
$$
\widetilde{q}(t,x):=1\wedge\frac{\rho(x)^{\alpha/2}}{\sqrt{t}}.
$$
The following easy result will be used several times below.

\bl
For every $-1<\gamma<d/\alpha$, $t>0$ and $x\in\mR^d$,
\begin{align} \label{gam}
\int_{0}^{t}\varrho^\gamma(s,x)\dif
s\preceq
\frac{1}{|x|^{d-\alpha\gamma}}\wedge \frac{t^{1+\gamma}}{|x|^{d+\alpha}}.
\end{align}
\el

\begin{proof}
If $|x|\geq t^{1/\alpha}$, we have
\begin{align*}
\int_{0}^{t}\varrho^\gamma(s,x)\dif s\leq \int_{0}^{t}\frac{s^{\gamma}}{|x|^{d+\alpha}}\dif s\preceq \frac{t^{1+\gamma}}{|x|^{d+\alpha}}.
\end{align*}
If $|x|< t^{1/\alpha}$, we have
\begin{align*}
\int_{0}^{t}\varrho^\gamma(s,x)\dif s\leq \int_{0}^{|x|^\alpha}\frac{s^{\gamma}}{|x|^{d+\alpha}}\dif s+\int_{|x|^\alpha}^{\infty}s^{\gamma-\frac{d+\alpha}{\alpha}}\dif s\preceq\frac{1}{|x|^{d-\alpha\gamma}}.
\end{align*}
Combining the above computations, we get the desired result.
\end{proof}

In the remainder of this section, we fix an arbitrary $T>0$ and
assume that $D$ is a $C^{1,\theta}$ open set in $\mR^d$ with
$\theta\in (\alpha/2, 1]$. Recall that $p^D(t,x,y)$ is the transition density of $X^D$
and it holds that (see \cite{KK}) for all $(t,x,y)\in(0,T]\times D\times D$,
\begin{align}
p^{D}(t,x,y)\asymp
\widetilde{q}(t,x)\widetilde{q}(t,y)\varrho^1(t,x-y).\label{pd}
\end{align}
Although the classical 3-P inequality of the form \eqref{class3p} does not hold for $p^D(t,x,y)$, we do have the following generalized 3-P type inequality.

\bl
For any $0<s\le t\le T$ and $x, y, z\in D$, we have
\begin{align}
\frac{p^{D}(t-s,x,z)p^{D}(s,z,y)}{p^{D}(t,x,y)}\preceq
\rho(z)^\alpha\Big(\varrho^0(t-s,x-z)+\varrho^0(s,z-y)\Big).
\label{3pd}
\end{align}
\el
\begin{proof}
Note that
\begin{align*}
\big(|x-y|+t^{1/\alpha}\big)^{d+\alpha}\preceq \big(|x-z|+(t-s)^{1/\alpha}\big)^{d+\alpha}+\big(|z-y|+s^{1/\alpha}\big)^{d+\alpha}.
\end{align*}
Thus
\begin{align}
\frac{\varrho^1(t-s,x-z)\varrho^1(s,z-y)}{\varrho^1(t,x-y)}&=\frac{(t-s) s}{t}\cdot\frac{\varrho^0(t-s,x-z)\varrho^0(s,z-y)}{\varrho^0(t,x-y)}\no\\
&\preceq
\frac{(t-s) s}{t}\cdot\Big(\varrho^0(t-s,x-z)+\varrho^0(s,z-y)\Big).\label{3p}
\end{align}
It is obvious that
$$
\sqrt{s}\cdot \widetilde q(s,x)\leq\sqrt{t}\cdot \widetilde q(t,x).
$$
Combining this with \eqref{pd} we can derive that
\begin{align*}
&\frac{p^{D}(t-s,x,z)p^{D}(s,z,y)}{p^{D}(t,x,y)}\preceq\frac{\rho(z)^{\alpha}}{\sqrt{(t-s)s}}\cdot\frac{t}{\sqrt{(t-s)s}}
\cdot \frac{\varrho^1(t-s,x-z)\varrho^1(s,z-y)}{\varrho^1(t,x-y)}\\
&\preceq \frac{\rho(z)^{\alpha}}{\sqrt{(t-s)s}}\cdot\frac{t}{\sqrt{(t-s)s}}\cdot\frac{(t-s)s}{t}\cdot\Big(\varrho^0(t-s,x-z)+\varrho^0(s,z-y)\Big)\\
&=\rho(z)^\alpha\Big(\varrho^0(t-s,x-z)+\varrho^0(s,z-y)\Big).
\end{align*}
The proof is complete.
\end{proof}

We will also need the following generalized integral
inequality.

\bl
For any $t\in (0, T]$ and $y,z\in D$, we have
\begin{align}\label{2.6}
\widetilde{q}(t,z)\int^{t/2}_{0}\!\!s^{-1/\alpha}p^{D}(s,z,y)\dif s
\preceq
\widetilde{q}(t,y)\int^{t/2}_{0}\!\!s^{-1/\alpha}\widetilde{q}(s,z)\varrho^1(s,z-y)\dif s.
\end{align}
\el
\begin{proof}
It can be easily checked that (\ref{2.6}) holds when $\rho(y)\geq (t/2)^{1/\alpha}$
or $\rho(z)\leq2\rho(y)$. So we will assume
$\rho(y)<(t/2)^{1/\alpha}\wedge(\rho(z)/2)$ throughout this proof. Note
that in this case, we have
\begin{align*}
|z-y|\geq\rho(z)-\rho(y)\geq\frac{\rho(z)}{2}\geq\rho(y).
\end{align*}
For convenience, we define
\begin{align*}
\mathbf{L}:=\widetilde{q}(t,z)\int^{t/2}_{0}\!\!s^{-1/\alpha}p^{D}(s,z,y)\dif
s
\end{align*}
and
\begin{align*}
\mathbf{R}:=\widetilde{q}(t,y)\int^{t/2}_{0}\!\!s^{-1/\alpha}\widetilde{q}(s,z)\varrho^1(s,z-y)\dif
s.
\end{align*}
We deal with  three different cases separately.\\
Case 1: $(t/2)^{1/\alpha}<\rho(z)/2<|z-y|$.
In this case, we have
\begin{align}
\mathbf{L}&\preceq\int_{0}^{t/2}s^{-1/\alpha}
\left(1\wedge\frac{\rho(y)^{\alpha/2}}{\sqrt{s}}\right)\frac{s}{|z-y|^{d+\alpha}}\dif
s\no\\
&=\int_{0}^{\rho(y)^{\alpha}}\!s^{-1/\alpha}\frac{s}{|z-y|^{d+\alpha}}\dif
s+\int_{\rho(y)^{\alpha}}^{t/2}s^{-1/\alpha}\frac{\rho(y)^{\alpha/2}}{\sqrt{s}}\cdot\frac{s}{|z-y|^{d+\alpha}}\dif
s\no\\
&\asymp\frac{\rho(y)^{\alpha/2}}{|z-y|^{d+\alpha}}t^{3/2-1/\alpha}\label{left}
\end{align}
and
\begin{align*}
\mathbf{R}\asymp\frac{\rho(y)^{\alpha/2}}{\sqrt{t}}\int_{0}^{t/2}s^{-1/\alpha}\frac{s}{|z-y|^{d+\alpha}}\dif
s\asymp\frac{\rho(y)^{\alpha/2}}{|z-y|^{d+\alpha}}t^{3/2-1/\alpha}.
\end{align*}
Thus, we have $\mathbf{L}\preceq\mathbf{R}$ in this case.\\
Case 2: $\rho(z)/2\leq (t/2)^{1/\alpha}<|z-y|$.
By using the same argument as in \eqref{left}, we can get that
\begin{align*}
\mathbf{L}&\preceq\frac{\rho(z)^{\alpha/2}}{\sqrt{t}}\int_{0}^{t/2}s^{-1/\alpha}\left(1\wedge\frac{\rho(y)^{\alpha/2}}{\sqrt{s}}\right)\frac{s}{|z-y|^{d+\alpha}}\dif
s\\
&\asymp\frac{\rho(z)^{\alpha/2}}{\sqrt{t}}\cdot\frac{\rho(y)^{\alpha/2}}{|z-y|^{d+\alpha}}t^{3/2-1/\alpha}
\end{align*}
and
\begin{align*}
\mathbf{R}&\asymp\frac{\rho(y)^{\alpha/2}}{\sqrt{t}}\int_{0}^{t/2}s^{-1/\alpha}\left(1\wedge\frac{\rho(z)^{\alpha/2}}{\sqrt{s}}\right)\frac{s}{|z-y|^{d+\alpha}}\dif
s\\
&\asymp\frac{\rho(z)^{\alpha/2}}{\sqrt{t}}\cdot\frac{\rho(y)^{\alpha/2}}{|z-y|^{d+\alpha}}t^{3/2-1/\alpha}.
\end{align*}
Thus, we also have $\mathbf{L}\preceq \mathbf{R}$ in this case.\\
Case 3: $\rho(z)/2\leq|z-y|\leq (t/2)^{1/\alpha}$.
In this case, we have
\begin{align*}
\mathbf{L}&\preceq\frac{\rho(z)^{\alpha/2}}{\sqrt{t}}\int_{0}^{t/2}s^{-1/\alpha}\left(1\wedge\frac{\rho(y)^{\alpha/2}}{\sqrt{s}}\right)\left(s^{-d/\alpha}\wedge\frac{s}{|z-y|^{d+\alpha}}\right)\dif
s\\
&\asymp\frac{\rho(z)^{\alpha/2}}{\sqrt{t}}\bigg(\int_{0}^{\rho(y)^{\alpha}}\!\!s^{-1/\alpha}\frac{s}{|z-y|^{d+\alpha}}\dif
s+\int_{\rho(y)^{\alpha}}^{|z-y|^{\alpha}}s^{-1/\alpha}\frac{\rho(y)^{\alpha/2}}{\sqrt{s}}\cdot\frac{s}{|z-y|^{d+\alpha}}\dif
s\\
&\quad\quad\quad\quad\quad\quad+\int_{|z-y|^{\alpha}}^{t/2}s^{-1/\alpha}\frac{\rho(y)^{\alpha/2}}{\sqrt{s}}s^{-d/\alpha}\dif
s\bigg)\\
&\asymp\frac{\rho(z)^{\alpha/2}\cdot\rho(y)^{\alpha/2}}{\sqrt{t}\cdot|z-y|^{d+1-\alpha/2}}-\frac{\rho(z)^{\alpha/2}\cdot\rho(y)^{\alpha/2}}{t^{(d+1)/\alpha}}.
\end{align*}
Using the same idea, we can also get
\begin{align*}
\mathbf{R}&\asymp\frac{\rho(y)^{\alpha/2}}{\sqrt{t}}\int_{0}^{t/2}s^{-1/\alpha}\left(1\wedge\frac{\rho(z)^{\alpha/2}}{\sqrt{s}}\right)\left(s^{-d/\alpha}\wedge\frac{s}{|z-y|^{d+\alpha}}\right)\dif
s\\
&\asymp\frac{\rho(z)^{\alpha/2}\cdot\rho(y)^{\alpha/2}}{\sqrt{t}\cdot|z-y|^{d+1-\alpha/2}}-\frac{\rho(z)^{\alpha/2}\cdot\rho(y)^{\alpha/2}}{t^{(d+1)/\alpha}}.
\end{align*}
Thus, $\mathbf{L}\preceq \mathbf{R}$ in true.
The proof is now complete.
\end{proof}

\section{Proof of Theorem \ref{main}}

Throughout this section, unless specified otherwise, we always assume that
$b: D\to \mR^d$ belongs to $\mK_D^{\alpha-1}$.
The following lemma plays an important role in proving our main result.

\bl\label{lemma3.2}
Let $T>0$.
For any $t\in (0, T]$, there exists a
constant $C(t)=C(t, b)>0$ such that
all $x, y\in D$, we have
\begin{align*}
\int_0^t\!\!\!\int_{D}p^{D}(t-s,x,z)|b(z)|\cdot|\nabla_z
p^{D}(s,z,y)|\dif z\dif s\leq
C(t)p^{D}(t,x,y),
\end{align*}
where $C(t)$ is nondecreasing in $t$ and $C(t)\rightarrow0$ as $t\rightarrow0$.
\el
\begin{proof}
Define
\begin{align*}
\cI&:= \frac{1}{p^{D}(t,x,y)}\int_0^tp^{D}(t-s,x,z)\cdot|\nabla_z p^{D}(s,z,y)|\dif s.
\end{align*}
Then, by (\ref{grad2-a}), we have
\begin{align*}
\cI&\preceq \int_0^t\frac{1}{\rho(z)}1_{[\rho(z)<s^{1/\alpha}\wedge(t-s)^{1/\alpha}]}\frac{p^{D}(t-s,x,z)p^{D}(s,z,y)}{p^{D}(t,x,y)}\dif s\\
&\qquad +\int_0^t\frac{1}{\rho(z)\wedge
s^{1/\alpha}}1_{[\rho(z)\geq s^{1/\alpha}\wedge(t-s)^{1/\alpha}]}\frac{p^{D}(t-s,x,z)p^{D}(s,z,y)}{p^{D}(t,x,y)}\dif
s\\
&=:\cI_1+\cI_2.
\end{align*}
On one hand, we have by \eqref{3pd} that
\begin{align*}
\cI_1
&\preceq\int_0^t\rho(z)^{\alpha-1}1_{[\rho(z)<s^{1/\alpha}\wedge(t-s)^{1/\alpha}]}\Big(\varrho^{0}(t-s,x-z)+\varrho^{0}(s,z-y)\Big)\dif s\\
&\leq\int_0^t\Big((t-s)^{1-1/\alpha}\varrho^{0}(t-s,x-z)\dif s+\int_0^ts^{1-1/\alpha}\varrho^{0}(s,z-y)\Big)\dif s\\
&=\int_{0}^{t}s^{-1/\alpha}\Big(\varrho^{1}(s,x-z)+\varrho^{1}(s,z-y)\Big)\dif s.
\end{align*}
We proceed to show that  $\cI_2$ has the same estimate:
\begin{align}
\cI_2\preceq \int_0^ts^{-1/\alpha}\Big(\varrho^1(s,x-z)+\varrho^1(s,z-y)\Big)\dif s.    \label{iii}
\end{align}
Since
$$
\frac{1}{\rho(z)\wedge
s^{1/\alpha}}1_{[\rho(z)\geq s^{1/\alpha}\wedge(t-s)^{1/\alpha}]}\leq\frac{1}{s^{1/\alpha}}+\frac{1}{(t-s)^{1/\alpha}},
$$
we have
\begin{align*}
\cI_2\preceq
\int_0^t\left(\frac{1}{s^{1/\alpha}}+\frac{1}{(t-s)^{1/\alpha}}\right)\frac{p^{D}(t-s,x,z)p^{D}(s,z,y)}{p^{D}(t,x,y)}\dif s.
\end{align*}
Using the symmetry in $s$ and $t-s$, we only need to prove that
\begin{align*}
\hat\cI_2:=\!\!\int_0^t\frac{1}{s^{1/\alpha}}\frac{p^{D}(t-s,x,z)p^{D}(s,z,y)}{p^{D}(t,x,y)}\dif
s\!
\preceq\!\!\int_{0}^{t}\!\frac{1}{s^{1/\alpha}}\Big(\varrho^1(s,x-z)+\varrho^1(s,z-y)\Big)\dif
s.
\end{align*}
By (\ref{2.6}), we have
\begin{align*}
\hat\cI_2&\preceq\frac{p^{D}(t,x,z)}{p^{D}(t,x,y)}\int_{0}^{t/2}\!s^{-1/\alpha}p^D(s,z,y)\dif s+\frac{p^{D}(t,z,y)}{p^{D}(t,x,y)}\int_{t/2}^{t}s^{-1/\alpha}p^D(t-s,x,z)\dif s\\
&\preceq\frac{\widetilde{q}(t,x)\widetilde{q}(t,y)}{p^{D}(t,x,y)}\varrho^1(t,x-z)\int_{0}^{t/2}\!s^{-1/\alpha}\widetilde{q}(s,z)\varrho^1(s,z-y)\dif s\\
&\qquad+\frac{\widetilde{q}(t,x)\widetilde{q}(t,y)}{p^{D}(t,x,y)}\varrho^1(t,z-y)\int_{0}^{t/2}\!s^{-1/\alpha}\widetilde{q}(s,z)\varrho^1(s,x-z)\dif
s\\
&=:\hat\cI_{21}+\hat\cI_{22},
\end{align*}
where we have used a change of variables and the facts that
$\rho^1(t, x-z)\asymp \rho^1(t-s, x-z)$ for $s\in (0, t/2)$.
It suffices to take care of one of the two terms of the right hand side
above, the other term can be handled in a similar fashion.
By \eqref{3p}, we have
\begin{align*}
\hat\cI_{21}
&\preceq\frac{\widetilde{q}(t,x)\widetilde{q}(t,x)}{p^{D}(t,x,y)}\int_{0}^{t/2}s^{-1/\alpha}\varrho^1(t-s,x-z)\varrho^1(s,z-y)\dif s\\
&\preceq\int_{0}^{t/2}s^{-1/\alpha}\frac{\varrho^1(t-s,x-z)\varrho^1(s,z-y)}{\varrho^1(t,x-y)}\dif
s\\
&\preceq\int_{0}^{t/2}s^{-1/\alpha}\left(\frac{s}{\left(|x-z|+(t-s)^{1/\alpha}\right)^{d+\alpha}}+\frac{s}{\left(|z-y|+s^{1/\alpha}\right)^{d+\alpha}}\right)\dif s\\
&\leq\int_{0}^{t/2}s^{-1/\alpha}\Big(\varrho^1(s,x-z)+\varrho^1(s,z-y)\Big)\dif
s.
\end{align*}
Combining the four displays above, we get \eqref{iii}. Using (\ref{gam}) with $\gamma=1-1/\alpha$, we arrive at
\begin{align}
\cI&\preceq \int_0^ts^{-1/\alpha}\Big(\varrho^1(s,x-z)+\varrho^1(s,z-y)\Big)\dif s\nonumber\\
&\preceq \left(\frac{1}{|x-z|^{d+1-\alpha}}\wedge\frac{t^{2-1/\alpha}}{|x-z|^{d+\alpha}}+\frac{1}{|z-y|^{d+1-\alpha}}\wedge\frac{t^{2-1/\alpha}}{|z-y|^{d+\alpha}}\right).
\label{e:I}\end{align}
Consequently,
\begin{align*}
&\int_0^t\!\!\!\int_{D}p^{D}(t-s,x,z)|b(z)|\cdot|\nabla_z
p^{D}(s,z,y)|\dif z\dif s\\
&\preceq
\sup_{w\in D}\int_D\left(\frac{1}{|w-z|^{d+1-\alpha}}\wedge
\frac{t^{2-1/\alpha}}{|w-z|^{d+\alpha}}\right)|b(z)|\dif z \cdot p^{D}(t,x,y).
\end{align*}
The desired conclusion now follows from Lemma \ref{include} with $\beta=2-1/\alpha$.
\end{proof}

To derive our gradient estimate, we will also need the following result.

\bl\label{lemma3.3}
Let $T>0$. For
any $t\in(0,T]$, there exists $\widehat{C}(t)=\widehat{C}(t,b)>0$  such that for all
$x, y\in D$, we have
\begin{align*}
\int_0^t\!\!\!\int_{D}|\nabla_x p^{D}(t-s,x,z)||b(z)|\cdot|\nabla_z
p^{D}(s,z,y)|\dif z\dif s\leq
\frac{\widehat{C}(t)}{\rho(x)\wedge
t^{1/\alpha}}p^{D}(t,x,y),
\end{align*}
where $\widehat{C}(t)$ is nondecreasing in $t$ and
$\widehat{C}(t)\rightarrow0$ as $t\rightarrow0$.
\el
\begin{proof}
Define
\begin{align*}
\cQ&:=\frac{\rho(x)\wedge
t^{1/\alpha}}{p^{D}(t,x,y)}\int_0^t|\nabla_x p^{D}(t-s,x,z)|\cdot|\nabla_z
p^{D}(s,z,y)|\dif s.
\end{align*}
By (\ref{grad2-a}), we have that
\begin{align*}
\cQ&\preceq\frac{\rho(x)\wedge
t^{1/\alpha}}{p^{D}(t,x,y)}\int_0^t\frac{p^{D}(t-s,x,z)}{\rho(x)\wedge(t-s)^{1/\alpha}}
\cdot |\nabla_z p^{D}(s,z,y)|\dif s\\
&\preceq\frac{\rho(x)\wedge
t^{1/\alpha}}{p^{D}(t,x,y)}\int_0^t\frac{p^{D}(t-s,x,z)}{\rho(x)}1_{[\rho(x)\leq(t-s)^{1/\alpha}]}\cdot
|\nabla_z p^{D}(s,z,y)|\dif s\\
&\qquad+\frac{\rho(x)\wedge
t^{1/\alpha}}{p^{D}(t,x,y)}\int_0^t\frac{p^{D}(t-s,x,z)}{(t-s)^{1/\alpha}}1_{[\rho(x)>(t-s)^{1/\alpha}]}\cdot\frac{p^{D}(s,z,y)}{\rho(z)\wedge s^{1/\alpha}}\dif s\\
&=:\cQ_{1}+\cQ_{2}.
\end{align*}
Using \eqref{e:I} in the second line below, we get
that
\begin{align*}
\cQ_{1}&\leq\frac{1}{p^{D}(t,x,y)}\int_0^tp^{D}(t-s,x,z)\cdot
|\nabla_z p^{D}(s,z,y)|\dif s\\
&\preceq \left(\frac{1}{|x-z|^{d+1-\alpha}}\wedge\frac{t^{2-1/\alpha}}{|x-z|^{d+\alpha}}+\frac{1}{|z-y|^{d+1-\alpha}}\wedge\frac{t^{2-1/\alpha}}{|z-y|^{d+\alpha}}\right).
\end{align*}
To deal with $\cQ_{2}$, we rewrite it as
\begin{align*}
\cQ_{2}&=\frac{\rho(x)\wedge
t^{1/\alpha}}{p^{D}(t,x,y)}\int_0^t\frac{p^{D}(t-s,x,z)}{(t-s)^{1/\alpha}}1_{[\rho(x)>(t-s)^{1/\alpha}]}\cdot\frac{p^{D}(s,z,y)}{s^{1/\alpha}}1_{[\rho(z)> s^{1/\alpha})]}\dif s\\
&\qquad+\frac{\rho(x)\wedge
t^{1/\alpha}}{p^{D}(t,x,y)}\int_0^t\frac{p^{D}(t-s,x,z)}{(t-s)^{1/\alpha}}1_{[\rho(x)>(t-s)^{1/\alpha}]}\cdot\frac{p^{D}(s,z,y)}{\rho(z)}1_{[\rho(z)\leq s^{1/\alpha})]}\dif s\\
&=:\cQ_{21}+\cQ_{22}.
\end{align*}
On one hand, we have by (\ref{2.6}) that
\begin{align*}
\cQ_{21}&\leq\frac{\rho(x)\wedge
t^{1/\alpha}}{p^{D}(t,x,y)}\left(\int_0^{t/2}+\int_{t/2}^{t}\right)\frac{p^{D}(t-s,x,z)}{(t-s)^{1/\alpha}}\cdot\frac{p^{D}(s,z,y)}{s^{1/\alpha}}\dif s\\
&\preceq\frac{\rho(x)\wedge
t^{1/\alpha}}{p^{D}(t,x,y)}\frac{\widetilde{q}(t,x)\widetilde{q}(t,y)}{t^{1/\alpha}}\varrho^1(t,x-z)\int_{0}^{t/2}s^{-1/\alpha}\tilde{q}(s,z)\varrho^1(s,z-y)\dif s\\
&\qquad+\frac{\rho(x)\wedge
t^{1/\alpha}}{p^{D}(t,x,y)}\frac{\widetilde{q}(t,x)\widetilde{q}(t,y)}{t^{1/\alpha}}\varrho^1(t,z-y)\int_{0}^{t/2}s^{-1/\alpha}\tilde{q}(s,z)\varrho^1(s,x-z)\dif s.
\end{align*}
Repeating the argument used to estimate $\hat\cI_{21}$ in the proof of Lemma \ref{lemma3.2},
we get that
\begin{align*}
\cQ_{21}&\leq\frac{\widetilde{q}(t,x)\widetilde{q}(t,y)}{p^{D}(t,x,y)}\int_{0}^{t/2}s^{-1/\alpha}\varrho^1(t-s,x-z)\varrho^1(s,z-y)\dif s\\
&\preceq\left(\frac{1}{|x-z|^{d+1-\alpha}}\wedge\frac{t^{2-1/\alpha}}{|x-z|^{d+\alpha}}+\frac{1}{|z-y|^{d+1-\alpha}}\wedge\frac{t^{2-1/\alpha}}{|z-y|^{d+\alpha}}\right).
\end{align*}
To deal with $\cQ_{22}$, we write
\begin{align*}
\cQ_{22}&=\frac{\rho(x)\wedge
t^{1/\alpha}}{p^{D}(t,x,y)}\left(\int_0^{t/2}+\int_{t/2}^{t}\right)\frac{p^{D}(t-s,x,z)}{(t-s)^{1/\alpha}}1_{[\rho(x)>(t-s)^{1/\alpha}]}
\cdot\frac{p^{D}(s,z,y)}{\rho(z)}1_{[\rho(z)\leq s^{1/\alpha})]}\dif s\\
&=:\hat{\cQ}_{21}+\hat{\cQ}_{22}.
\end{align*}
We can use \eqref{3pd} to deduce that
\begin{align*}
\hat{\cQ}_{21}&\leq\frac{\rho(x)\wedge t^{1/\alpha}}{t^{1/\alpha}}\int_{0}^{t/2}\rho(z)^{\alpha-1}1_{[\rho(z)\leq s^{1/\alpha}]}\Big(\varrho^0(t-s,x-z)+\varrho^0(s,z-y)\Big)\dif s\\
&\leq\int_{0}^{t/2}s^{1-1/\alpha}\Big(\varrho^0(s,x-z)+\varrho^0(s,z-y)\Big)\dif
s\\
&\preceq\left(\frac{1}{|x-z|^{d+1-\alpha}}\wedge\frac{t^{2-1/\alpha}}{|x-z|^{d+\alpha}}+\frac{1}{|z-y|^{d+1-\alpha}}\wedge\frac{t^{2-1/\alpha}}{|z-y|^{d+\alpha}}\right).
\end{align*}
We claim that
\begin{align}
\hat\cQ_{22}\preceq\left(\frac{1}{|x-z|^{d+1-\alpha}}\wedge\frac{t^{2-2/\alpha}}{|x-z|^{d+\alpha-1}}+\frac{1}{|z-y|^{d+1-\alpha}}\wedge\frac{t^{2-2/\alpha}}{|z-y|^{d+\alpha-1}}\right). \label{cla}
\end{align}
To prove this claim, we write
\begin{align*}
\hat{\cQ}_{22}&=\frac{\rho(x)\wedge
t^{1/\alpha}}{p^{D}(t,x,y)}\int_{t/2}^{t}\frac{p^{D}(t-s,x,z)}{(t-s)^{1/\alpha}}1_{[\rho(x)>(t-s)^{1/\alpha}]}\cdot\frac{p^{D}(s,z,y)}{\rho(z)}
1_{[(t-s)^{1/\alpha}<\rho(z)\leq s^{1/\alpha})]}\dif s\\
&\qquad+\frac{\rho(x)\wedge
t^{1/\alpha}}{p^{D}(t,x,y)}\int_{t/2}^{t}\frac{p^{D}(t-s,x,z)}{(t-s)^{1/\alpha}}1_{[\rho(x)>(t-s)^{1/\alpha}]}\cdot\frac{p^{D}(s,z,y)}{\rho(z)}
1_{[\rho(z)\leq (t-s)^{1/\alpha})]}\dif s\\
&=:\tilde{\sQ}_{1}+\tilde{\sQ}_{2}.
\end{align*}
If we denote $A:=[(t-s)^{1/\alpha}<\rho(z)\leq s^{1/\alpha}]$, then
we have by (\ref{3p}) that
\begin{align*}
\tilde{\sQ}_{1}&\preceq\frac{\rho(x)\wedge
t^{1/\alpha}}{\sqrt{t}\cdot\tilde{q}(t,x)}\int_{t/2}^{t}(t-s)^{1-1/\alpha}\rho(z)^{\alpha/2-1}\Big(\varrho^{0}(t-s,x-z)+\varrho^{0}(s,z-y)\Big)1_{A}\dif s\\
&\leq
t^{1/\alpha-1/2}\int_{t/2}^{t}(t-s)^{1-1/\alpha}\rho(z)^{\alpha/2-1}\varrho^{0}(s,z-y)1_{A}\dif
s\\
&\qquad+\rho(x)^{1-\alpha/2}\int_{t/2}^{t}(t-s)^{1-1/\alpha}\rho(z)^{\alpha/2-1}\varrho^{0}(t-s,x-z)1_{A}\dif
s\\
&=:\tilde{\sQ}_{11}+\tilde{\sQ}_{12},
\end{align*}
where in the second inequality we used the fact that $\frac{\rho(x)\wedge
t^{1/\alpha}}{\widetilde{q}(t,x)}=\sqrt{t}\big(\rho(x)\wedge
t^{1/\alpha}\big)^{1-\alpha/2}$. One can easily check  that
\begin{align*}
\tilde{\sQ}_{11}
&\preceq
t^{1/\alpha-1/2}\varrho^{0}(t,z-y)\int_{t/2}^{t}(t-s)^{3/2-2/\alpha}\dif
s\\
&\preceq\frac{1}{|z-y|^{d+1-\alpha}}\wedge\frac{t}{|z-y|^{d+1}}.
\end{align*}
By the fact $\rho(x)\leq\rho(z)+|x-z|,$ we further have
\begin{align*}
\tilde{\sQ}_{12}&\preceq\big(\rho(z)^{1-\alpha/2}+|x-z|^{1-\alpha/2}\big)\int_{t/2}^{t}(t-s)^{-1/\alpha}\rho(z)^{\alpha/2-1}\varrho^{1}(t-s,x-z)1_{A}\dif
s\\
&\leq\int_{t/2}^{t}(t-s)^{-1/\alpha}\varrho^{1}(t-s,x-z)\dif
s+|x-z|^{1-\alpha/2}\!\!\!\int_{t/2}^{t}(t-s)^{1/2-2/\alpha}\varrho^{1}(t-s,x-z)\dif s\\
&\preceq\frac{1}{|x-z|^{d+1-\alpha}}\wedge
\frac{t^{2-2/\alpha}}{|x-z|^{d+\alpha-1}}.
\end{align*}
To estimate $\tilde\sQ_2$, we can use (\ref{3pd}) to deduce that
\begin{align*}
\tilde{\sQ}_{2}&\preceq\big(\rho(x)\wedge
t^{1/\alpha}\big)\!\!\int_{t/2}^{t}(t-s)^{-1/\alpha}\rho(z)^{\alpha-1}1_{[\rho(z)\leq(t-s)^{1/\alpha}]}\Big(\varrho^{0}(t-s,x-z)+\varrho^{0}(s,z-y)\Big)\dif
s\\
&\leq
t^{1/\alpha}\!\!\int_{t/2}^{t}(t-s)^{-1/\alpha}\rho(z)^{\alpha-1}1_{[\rho(z)\leq(t-s)^{1/\alpha}]}\varrho^{0}(s,z-y)\dif
s\\
&\qquad+\rho(x)\!\int_{t/2}^{t}(t-s)^{-1/\alpha}\rho(z)^{\alpha-1}1_{[\rho(z)\leq(t-s)^{1/\alpha}]}\varrho^{0}(t-s,x-z)\dif
s.
\end{align*}
Then by the same argument used for $\tilde{\sQ}_{1}$,
one can check  that
\begin{align*}
\tilde{\cQ}_{2}\preceq\frac{1}{|z-y|^{d+1-\alpha}}\wedge\frac{t}{|z-y|^{d+1}}+\frac{1}{|x-z|^{d+1-\alpha}}\wedge\frac{t^{2-2/\alpha}}{|x-z|^{d+\alpha-1}}.
\end{align*}
Since for any $x,z\in D$, we have
\begin{align}\label{e:neweq}
\frac{1}{|x-z|^{d+1-\alpha}}\wedge\frac{t^{2-1/\alpha}}{|x-z|^{d+\alpha}}
&\leq\frac{1}{|x-z|^{d+1-\alpha}}\wedge\frac{t}{|x-z|^{d+1}}\no\\
&\leq\frac{1}{|x-z|^{d+1-\alpha}}\wedge\frac{t^{2-2/\alpha}}{|x-z|^{d+\alpha-1}},
\end{align}
combining the displays above, we get (\ref{cla}).

Combining \eqref{e:neweq}, \eqref{cla} with our estimates for $\cQ_1$, $\cQ_{21}$,
$\hat{\cQ}_{21}$,  we get
$$
\cQ\preceq\left(\frac{1}{|x-z|^{d+1-\alpha}}\wedge\frac{t^{2-2/\alpha}}{|x-z|^{d+\alpha-1}}+\frac{1}{|z-y|^{d+1-\alpha}}\wedge\frac{t^{2-2/\alpha}}{|z-y|^{d+\alpha-1}}\right).
$$
Hence
\begin{align*}
&\int_0^t\!\!\!\int_{D}|\nabla_x p^{D}(t-s,x,z)||b(z)|\cdot|\nabla_z
p^{D}(s,z,y)|\dif z\dif s\\
&\preceq\sup_{w\in D}\int_D\left(\frac{1}{|w-z|^{d+1-\alpha}}\wedge\frac{t^{2-2/\alpha}}{|w-z|^{d+\alpha-1}}\right)|b(z)|\dif z\cdot \frac{1}{\rho(x)\wedge
t^{1/\alpha}}p^{D}(t,x,y).
\end{align*}
The desired conclusion now follows from Lemma \ref{include} with $\beta=2-2/\alpha$.
\end{proof}

We now proceed to solve the integral equation (\ref{du}). For
all $(t,x,y)\in(0,T]\times D\times D$, set
$p_0(t,x,y):=p^D(t,x,y)$, and define inductively that for
$k\geq1$,
\begin{align}\label{induction}
p_k(t,x,y):=\int_0^t\!\!\!\int_Dp_{k-1}(t-s,x,z)b(z)\cdot\nabla_z
p_{0}(s,z,y)\dif z\dif s.
\end{align}
The following result is an easy consequence of Lemmas \ref{lemma3.2} and  \ref{lemma3.3}.

\bl
Let $T>0$. For
every $k\geq 1$ and $x,y\in D$, we have
\begin{align}\label{3.5}
|p_k(t,x,y)|\leq C(t)^{k}p^D(t,x,y)
\end{align}
and
\begin{align}\label{GK}
|\nabla_x p_k(t,x,y)|\leq \frac{\widehat{C}(t)^{k}}{\rho(x)\wedge
t^{1/\alpha}}p^{D}(t,x,y),
\end{align}
where $C(t)$ is the constant in Lemma \ref{lemma3.2} and $\widehat
C(t)$ is the constant in Lemma \ref{lemma3.3}. Moreover, it holds
that
\begin{align}
p_k(t,x,y)=\int_0^t\!\!\!\int_Dp_{0}(t-s,x,z)b(z)\cdot\nabla_z p_{k-1}(s,z,y)\dif z\dif s.    \label{pkpk}
\end{align}
\el
\begin{proof}
We first prove (\ref{3.5}) by induction. By Lemma \ref{lemma3.2}, we know that
(\ref{3.5}) holds for $k=1$. Now suppose that it holds for
$k>1$. Then by definition and using Lemma \ref{lemma3.2} again, we
have
\begin{align*}
|p_{k+1}(t,x,y)|
&\leq\int_0^t\!\!\!\int_D|p_{k}(t-s,x,z)|\cdot|b(z)|\cdot|\nabla_z
p_{0}(s,z,y)|\dif z\dif s\\
&\leq
{C(t)}^{k}\int_0^t\!\!\!\int_Dp^{D}(t-s,x,z)|b(z)|\cdot|\nabla_z
p_{0}(s,z,y)|\dif z\dif s\\
&\leq {C(t)}^{k+1}p^{D}(t,x,y).
\end{align*}
Following the same argument with Lemma \ref{lemma3.2} replaced by Lemma \ref{lemma3.3}, we can
show (\ref{GK}) is true. We proceed to prove (\ref{pkpk}). It is obvious that
(\ref{pkpk}) holds for $k=1$. Suppose that it is true for certain
$k>1$. Then, we have by Fubini's theorem that
\begin{align*}
&p_{k+1}(t,x,y)=\int_0^t\!\!\!\int_Dp_{k}(t-s,x,z)b(z)\cdot\nabla_z p_{0}(s,z,y)\dif z\dif s\\
&=\int_0^t\!\!\!\int_D\int_0^{t-s}\!\!\!\int_Dp_{0}(t-s-r,x,u)b(u)\cdot\nabla_u p_{k-1}(r,u,z)\dif u\dif r\\
&\qquad\qquad\qquad\qquad\qquad\qquad\qquad\times b(z)\cdot\nabla_z p_{0}(s,z,y)\dif z\dif s\\
&=\int_0^t\!\!\!\int_Dp_{0}(t-\hat{r},x,u)b(u)\!\cdot\!\int_0^{\hat{r}}\!\!\!\int_D\nabla_u p_{k-1}(\hat{r}-s,u,z)\!\cdot\! b(z)\!\cdot\!\nabla_z p_{0}(s,z,y)\dif z\dif s\dif u\dif \hat{r}\\
&=\int_0^t\!\!\!\int_Dp_{0}(t-\hat{r},x,u)b(u)\cdot\nabla_u
p_{k}(\hat{r},u,y)\dif u\dif \hat{r},
\end{align*}
here in the third equality, we used the change of variable $\hat{r}=r+s$. The proof is complete.
\end{proof}

Now, we are in the position to give:

\begin{proof}[Proof of Theorem \ref{main}]
Let $p_k$ be defined by \eqref{induction}.
It follows from Lemma \ref{lemma3.2} that there exists $T_0\in (0, 1]$ such that
$C(T_0)<1/4$. Hence
\begin{align}
\sum_{k=0}^{\infty}|p_{k}(t,x,y)|\leq
\frac43 p^{D}(t,x,y),
\quad(t,x,y)\in(0, T_0]\times D\times D,   \label{es}
\end{align}
which means that the series $\sum_{k=0}^{\infty}p_{k}(t,x,y)$ is
convergent on $(0, T_0]\times D\times D$. Define
$p^{b,D}(t,x,y):=\sum_{k=0}^{\infty}p_{k}(t,x,y)$ on $(0, T_0]\times
D\times D$. By \eqref{induction}, we have
\begin{equation}\label{solution}
\sum_{k=0}^{n+1}p_k(t,x,y)=p_0(t,x,y)+\int_0^t\!\!\!\int_D\sum_{k=0}^np_k(t-s,x,z)b(z)\cdot\nabla_z p_0(s,z,y)\dif z\dif s.
\end{equation}
Letting $n\rightarrow\infty$ on both sides, we get \eqref{du}.

\vspace{2mm} \noindent(i) The upper bound on $(0, T_0]\times D\times
D$ follows from \eqref{es}. As for the lower bound on $(0,
T_0]\times D\times D$, we have
\begin{align*}
p^{b,D}(t,x,y)\geq p^D(t,x,y)-\sum_{k=1}^\infty
|p_k(t,x,y)|\geq \frac23 p^D(t,x,y).
\end{align*}
Thus, \eqref{1.5} is valid on $(0, T_0]\times D\times D$.

Now let $\widetilde{p}^{b,D}(t,x,y)$ be another solution to
\eqref{du} satisfying \eqref{1.5}, with $T$ replaced by $T_0$.
We claim that for every $k\in \mathbb{N}$ and
$t\in(0,T_0], x,y\in D$, there exists a constant $C_0$ such that
\begin{align}\label{unique}
|p^{b,D}(t,x,y)-\widetilde{p}^{b,D}(t,x,y)|\leq C_0{C(t)}^{k}p^{D}(t,x,y).
\end{align}
Indeed, for $k=1$, using \eqref{du}, \eqref{1.5} and Lemma \ref{lemma3.2} we
have
\begin{align*}
&\quad|p^{b,D}(t,x,y)-\widetilde{p}^{b,D}(t,x,y)|\\
&\leq\int_{0}^{t}\!\!\!\int_{D}\big(|p^{b,D}(t-s,x,z)|+|\widetilde{p}^{b,D}(t-s,x,z)|\big)\cdot|b(z)|\cdot|\nabla_z
p^{D}(s,z,y)|\dif z\dif s\\
&\leq C_0\int_{0}^{t}\!\!\!\int_{D}p^{D}(t-s,x,z)\cdot|b(z)|\cdot|\nabla_z
p^{D}(s,z,y)|\dif z\dif s\leq C_0C(t)p^{D}(t,x,y).
\end{align*}
Suppose that \eqref{unique} holds for some $k\in \mathbb{N}$. By
\eqref{du}, Lemma \ref{lemma3.2} and the induction hypothesis, we
have
\begin{align*}
&\quad|p^{b,D}(t,x,y)-\widetilde{p}^{b,D}(t,x,y)|\\
&\leq\int_{0}^{t}\!\!\!\int_{D}|p^{b,D}(t-s,x,z)-\widetilde{p}^{b,D}(t-s,x,z)|\cdot|b(z)|\cdot|\nabla_z
p^{D}(s,z,y)|\dif z\dif s\\
&\leq C_0{C(t)}^{k}\int_{0}^{t}\!\!\!\int_{D}p^{D}(t-s,x,z)\cdot|b(z)|\cdot|\nabla_z
p^{D}(s,z,y)|\dif z\dif s\leq C_0{C(t)}^{k+1}p^{D}(t,x,y).
\end{align*}
Since $C(t)<1$, letting $k\rightarrow\infty$, we obtain the
uniqueness.

\vspace{2mm} \noindent(ii)
By choosing $T_0$ smaller if necessary, we can assume that $\widehat C(T_0)<1$.
It then follows from \eqref{GK} that for every $t\in(0,T_0]$ and
$x,y\in D$,
\begin{align*}
\left|\sum_{k=0}^{\infty}\nabla_x
p_{k}(t,x,y)\right|\preceq\frac{1}{\rho(x)\wedge
t^{1/\alpha}}p^{D}(t,x,y),
\end{align*}
which means that (\ref{vi}) is true. Moreover, by (\ref{pkpk}) and
Fubini's theorem, we have
\begin{align*}
p^{b,D}(t,x,y)\!=\!\sum_{k=0}^{\infty}p_{k}(t,x,y)&=p^{D}(t,x,y)+\sum_{k=0}^{\infty}\int_0^t\!\!\!\int_Dp_{0}(t-s,x,z)b(z)\cdot\nabla_z
p_{k}(s,z,y)\dif z\dif s\\
&=p^{D}(t,x,y)+\int_0^t\!\!\!\int_Dp_{0}(t-s,x,z)b(z)\cdot\nabla_z
p^{b,D}(s,z,y)\dif z\dif s,
\end{align*}
that is (\ref{Duhamel}).

\vspace{2mm} \noindent(iii) By Fubini's theorem, we have for all
$0<s<t\le T_0$,
$$
\int_Dp^{b,D}(t-s,x,z)p^{b,D}(s,z,y)\dif z=\sum_{n=0}^{\infty}\sum_{m=0}^n\int_Dp_m(t-s,x,z)p_{n-m}(s,z,y)\dif z.
$$
Thus, to prove \eqref{eq21} for $0<s<t\le T_0$, it suffices to show that for each $n\in\mN_0$,
\begin{align}
\sum_{m=0}^n\int_Dp_{m}(t-s,x,z)p_{n-m}(s,z,y)\dif z=p_n(t,x,y).  \label{ck}
\end{align}
It is clear that the above equality holds for $n=0$. Suppose now that it holds for some $n\in\mN$. Write
$$
\sum_{m=0}^{n+1}\int_Dp_{m}(t-s,x,z)p_{n+1-m}(s,z,y)\dif z=\cJ_1+\cJ_2,
$$
where
$$
\cJ_1:=\int_Dp_{n+1}(t-s,x,z)p_{0}(s,z,y)\dif z
$$
and
$$
\cJ_2:=\sum_{m=0}^{n}\int_Dp_{m}(t-s,x,z)p_{n+1-m}(s,z,y)\dif z.
$$
By \eqref{induction} and Fubini's theorem, we have
\begin{align*}
\cJ_{1}&=\int_{D}\left(\int_{0}^{t-s}\!\!\!\int_{D}p_{n}(t-s-r,x,u)b(u)\cdot
\nabla_u p_{0}(r,u,z)\dif u\dif r\right)p_{0}(s,z,y)\dif z\\
&=\int_{0}^{t-s}\!\!\!\int_{D}p_{n}(t-s-r,x,u)b(u)\cdot\left(\int_{D}\nabla_u
p_{0}(r,u,z)p_{0}(s,z,y)\dif z\right)\dif u\dif r\\
&=\int_{s}^{t}\!\!\!\int_{D}p_{n}(t-r,x,u)b(u)\cdot\nabla_u
p_{0}(\hat{r},u,y)\dif u\dif r.
\end{align*}
Similarly, by \eqref{induction} and the induction hypothesis, we
have
\begin{align*}
\cJ_{2}=\int_{0}^{s}\!\!\!\int_{D}p_{n}(t-r,x,u)b(u)\cdot\nabla_u
p_{0}(r,u,y)\dif u\dif r.
\end{align*}
Hence,
\begin{align*}
\cJ_{1}+\cJ_{2}=\int_{0}^{t}\!\!\!\int_{D}p_{n}(t-r,x,u)b(u)\cdot\nabla_u
p_{0}(r,u,y)\dif u\dif r=p_{n+1}(t,x,y),
\end{align*}
which gives (\ref{ck}).

\vspace{2mm}
We now extend the function $p^{b, D}(t, x, y)$ from $(0, T]\times D\times D$ to
$(0, \infty)\times D\times D$ via the Chapman-Kolmogorov equation. Then it is routine to
extend the already assertions (i), (ii) (iii) on $(0, T_0]$ to $(0, T]$ for any $T>0$.

\vspace{2mm}
\noindent(iv)
Let $P_t^Df(x):=\int_Dp^{D}(t,x,y)f(y)\dif y$. By \eqref{du}, we have for
all $f\in C_c^2(D)$, $t>0$ and $x\in D$,
\begin{align}\label{RT5}
P_{t}^{b,D}f(x)=P^{D}_{t}f(x)+\int^t_0\!
P_{t-s}^{b,D}(b\cdot\nabla P^{D}_{s}f)(x)\dif s.
\end{align}
It then follows that for all $f\in C^2_c(D)$, $t>0$ and
$x\in D$,
\begin{align}
P_{t}^{b,D}f(x)-f(x)&=P^{D}_{t}f(x)-f(x)+
\int^t_0\! P_{t-s}^{b,D}(b\cdot\nabla P^{D}_{s}f)(x)\dif s\no\\
&=\int_0^t\!P^{D}_{t-s}(\Delta^{\alpha/2}|_Df)(x)\dif s+\int^t_0\!
P_{t-s}^{b,D}(b\cdot\nabla P^{D}_{s}f)(x)\dif s.\label{RT4}
\end{align}
Using \eqref{RT5} and Fubini's theorem, we get that for all $f\in C^2_c(D)$,
$t>0$ and $x\in D$,
\begin{align*}
&\quad\int_{0}^{t}P^{b,D}_{t-s}(\Delta^{\alpha/2}|_{D}f)(x)\dif
s-\int_{0}^{t}P^{D}_{t-s}(\Delta^{\alpha/2}|_{D}f)(x)\dif s\\
&=\int_{0}^{t}\!\!\!\int_{0}^{t-s}P^{b,D}_{t-s-r}(b\cdot\nabla
P^{D}_{r}\Delta^{\alpha/2}|_{D}f)(x)\dif r\dif s\\
&=\int_{0}^{t}\!\!\!\int_{s}^{t}P^{b,D}_{t-\hat{r}}(b\cdot\nabla
P^{D}_{\hat{r}-s}\Delta^{\alpha/2}|_{D}f)(x)\dif \hat{r}\dif s\\
&=\int_{0}^{t}P^{b,D}_{t-\hat{r}}\left(b\cdot\nabla\int_{0}^{\hat{r}}P^{D}_{\hat{r}-s}(\Delta^{\alpha/2}|_{D}f)ds\right)(x)\dif
s\dif \hat{r}\\
&=\int_{0}^{t}P^{b,D}_{t-\hat{r}}\left(b\cdot\nabla(P^{D}_{\hat{r}}f-f)(x)\right)\dif
\hat{r}.
\end{align*}
Combining this with \eqref{RT4}, we obtain that for all $f\in C^2_c(D)$,
$t>0$ and $x\in D$,
\begin{align*}
P_{t}^{b,D}f(x)-f(x)&
=\int_{0}^{t}P^{b,D}_{t-s}\left(\Delta^{\alpha/2}|_{D}+b\cdot\nabla\right)f(x)\dif s
\end{align*}
which gives (\ref{eq23}) for all $f\in C^2_c(D)$,
$t>0$ and $x\in D$.

\vspace{2mm} \noindent(v)  Set $b_n(x):=[(-n)\vee (b(x)\wedge
n)]1_{D}(x)$, then $b_n\in \mK^{\alpha-1}$. Let $p^{b_n,D}$ be the
transition density of $X^{b_n, D}$. Then we have
\begin{equation*}
\int_{D}p^{b_n,D}(t,x,y)\dif y\leq1.
\end{equation*}
It follows from  Lemma \ref{dens} that $p^{b_n,D}$  satisfies the Duhamel formula \eqref{du}
with $b$ replaced by $b_n$, thus by (\ref{solution}) and the uniqueness of solution to the integral equation (\ref{du}), $p^{b_n, D}$ can be defined as $\sum^\infty_{k=0}p_{k, n}(t, x, y)$, where $p_{k, n}$ is defined via \eqref{induction} with $b$ replaced by $b_n$,
Since $|b_n|\leq |b|$ and $b_n\rightarrow b$ as $n\rightarrow\infty$, as in the proof of \cite[Theorem 4.1]{C-H-X-Z}, we also have that
$p^{b_n,D}(t,x,y)\rightarrow p^{b,D}(t,x,y)$, which in turn implies (\ref{conser}).

\vspace{2mm} \noindent(vi) Since $p^{D}(t,x,y)$ is the transition
density of $X^{D}$, for any uniformly continuous function $f(x)$
with compact supports, we have
\begin{align*}
\lim_{t\downarrow0}\|P_t^Df-f\|_\infty=0.
\end{align*}
Meanwhile, by (\ref{1.5}) and Lemma \ref{lemma3.2} we have
\begin{align*}
&\quad\left|\!\int_{D}\!\!\!\left(\int_{0}^{t}\!\!\!\int_{D}p^{b,D}(t-s,x,z)b(z)\cdot\nabla_z
p^{D}(s,z,y)\dif z\dif s\!\right)\!f(y)\dif y\right|\\
&\preceq \|f\|_{\infty}\int_{D}\!\!\!\left(\int_{0}^{t}\!\!\!\int_{D}p^{D}(t-s,x,z)|b(z)|\cdot|\nabla_z
p^{D}(s,z,y)|\dif z\dif s\!\right)\dif y\\
&\leq C(t)\|f\|_{\infty}\int_{D}p^{D}(t,x,y)\dif y\leq
C(t)\|f\|_{\infty},
\end{align*}
where $C(t)\rightarrow0$ as $t\rightarrow0$, which yields
(\ref{con}) by (\ref{du}). The whole proof is finished.
\end{proof}

Finally, following the idea in \cite{L-L-W} we can give:

\begin{proof}[Proof of Corollary \ref{cor}]
By the two-sided heat kernel estimates (\ref{1.5}), there exists a constant $C>0$ such that for every $t\in(0,1]$ and $x,y\in D$,
$$
\frac{p^{b,D}(t,x,z)}{p^{b,D}(t,y,z)}\leq C_0\frac{\widetilde{q}(t,x)}{\widetilde{q}(t,y)}\frac{\varrho^1(t,x-z)}{\varrho^1(t,y-z)}\leq C\left(1\vee\frac{\rho(x)}{\rho(y)}\right)^{\alpha/2}\left(1\vee\frac{|x-y|}{t^{1/\alpha}}\right)^{d+\alpha}.
$$
Therefore, for any non-negative function $f\in\sB_b(D)$, $t\in(0,1)$ and $x,y\in D$, we have
\begin{align*}
P_t^{b,D}f(x)&=\int_D\frac{p^{b,D}(t,x,z)}{p^{b,D}(t,y,z)}p^{b,D}(t,y,z)f(z)\dif z\\
&\leq \left(\sup_{z\in D}\frac{p^{b,D}(t,x,z)}{p^{b,D}(t,y,z)}\right)\int_Dp^{b,D}(t,y,z)f(z)\dif z\\
&\leq C \left(1\vee\frac{\rho(x)}{\rho(y)}\right)^{\alpha/2}\left(1\vee\frac{|x-y|}{t^{1/\alpha}}\right)^{d+\alpha}P_t^{b,D}f(y),
\end{align*}
thus (\ref{ha}) holds for $t\in(0,1]$. For $T>1$, we can write by (\ref{ck}) that
$$
P_T^{b,D}f(x)=P_1^{b,D}P_{T-1}^{b,D}f(x).
$$
This together with the above inequality yields the desired result.
\end{proof}

\end{document}